\documentclass[12pt,a4paper]{amsart}

\usepackage{amsmath,amssymb,amsfonts,amsthm,enumerate,color}

\newtheorem{lemma}{Lemma}
\newtheorem{remark}{Remark}

\theoremstyle{remark}
\theoremstyle{plain}\newtheorem{theorem}{Theorem}
\newtheorem{corollary}[theorem]{Corollary}
\newtheorem{prop}[theorem]{Proposition}

\title[A random graph with duplications]{Asymptotic properties of a
  random graph with duplications}
\author{\'Agnes Backhausz}
\address{Department of Probability Theory and Statistics\\
E\"otv\"os Lor\'and University\\P\'azm\'any P.~s.\ 1/C, H-1117
  Budapest, Hungary}  
 \email{agnes@math.elte.hu}
 \author{Tam\'as F.~M\'ori}
 \address{Department of Probability Theory and Statistics\\
 E\"otv\"os Lor\'and University\\P\'azm\'any P.~s. 1/C, 
 H-1117 Budapest, Hungary}
 \email{mori@math.elte.hu}
 \dedicatory{\upshape {\sc Department of Probability Theory and
  Statistics}, \\ {\sc E\"otv\"os Lor\'and University}\\ P\'azm\'any
  P.~s. 1/C, H-1117 Budapest, Hungary\\                        
 \textit{E-mail address:} \texttt{agnes@math.elte.hu,
  mori@math.elte.hu} }
\thanks{Supported by the Hungarian Scientific Research Fund -- OTKA K 108615}


\date{30 August 2013}

\begin{document}

\begin{abstract}
We deal with a random graph model  evolving in discrete time steps
by duplicating and deleting the edges of randomly chosen
vertices. We prove the existence of an a.s. asymptotic degree
distribution, with stretched exponential decay; more precisely, the
proportion of vertices of degree $d$ tends to some positive number
$c_d>0$ almost surely as the number of steps goes to infinity, and
$c_d\sim (e\pi)^{1/2} d^{1/4} e^{-2\sqrt d}$ holds as $d\to\infty$.  
\end{abstract}

%\begin{titlepage}
\maketitle
\thispagestyle{empty}
%\end{titlepage}


\section{Introduction}

In the last decades, inspired by the examination of large real
networks, various types of  random graph models with preferential
attachment dynamics (meaning that vertices with larger degree have
larger chance to get new edges as the graph evolves randomly) were
introduced and analysed. After some early work \cite{yule, sim, szym},
this started with the seminal papers of Barab\'asi and Albert
\cite[1999]{ba} and Bollob\'as, Riordan, Spencer and Tusn\'ady
\cite[2001]{bb}. Among many others, we may mention the model of Cooper
and Frieze for the Internet \cite[2003]{cf} or that of Sridharan, Yong
Gao, Kui Wu and Nastos \cite[2011]{sridharan} for social networks.  

An important feature of these graph sequences is the scale-free
property: the proportion of vertices of degree $d$ tends  
to some positive number $c_d$ almost surely as the number of steps
goes to infinity, and $c_d\sim Kd^{-\gamma}$ holds as $d\rightarrow
\infty$ (throughout this paper, $a_d\sim b_d$ means that
$a_d/b_d\rightarrow 1$ as $d\rightarrow \infty$). To put it in another
way, the asymptotic degree distribution $(c_d)$ is polynomially
decaying. See also \cite{ba, fal, notices} and the  
references therein about the scale-free property of the internet. 

However, scale-free property captures only the behavior of the degrees of vertices,
and does not examine other kinds of structures. For example, especially
in biological networks, e.g., proteomes, it happens that we can find
groups of vertices having a similar neighborhood, that is, most of
their neighbors are the same. One can say that these networks are
highly clustered; loosely speaking, there are large cliques, in which
almost every vertex is connected to almost every other one, and there
are only a few edges going between cliques.  

A simple way to generate cliques is duplication: when a new vertex is
added, we choose an old vertex randomly, and connect the new vertex to
the neighbors of the old one. In other words, the new vertex becomes
a copy of the old vertex. Note that if the old vertex is chosen
uniformly at random, then the probability that a vertex of degree $d$
gets a new edge is just the probability that one of its neighbors is
chosen, which is proportional to its actual degree. Hence this model
is also driven by a kind of preferential attachment dynamics.  

After the duplication, we can add some extra edges randomly, or we can
delete some of them to guarantee that the network remains sparse.  The
graph may still have some large cliques due to the duplication. 

Duplication is not only a technical step that proved to be useful: it
is inherent. ``This may be because duplication of the information in
the genome is a dominant evolutionary force in shaping biological
networks (like gene regulatory networks and protein--protein
interaction networks)'' \cite{chung}.  

These kinds of models -- where the duplicated vertex is chosen
uniformly at random -- were examined for example by Kim, Krapivsky,
Kahng, and Redner \cite[2002]{kim}. In their model the new vertex is
connected to each neighbor of the chosen one with probability
$1-\delta$, independently. In addition, the new vertex is connected to
each old one independently with probability $\beta/n$ at the $n$th
step ($\delta, \beta$ are the parameters of the model). Scale free
property is claimed for this model. However, Pastor-Satorras, Smith
and Sol\'e \cite[2003]{pastor} stated that, instead of polynomial
decay, for the limit $c_d$ of the expected value of the proportion of
vertices of degree $d$ we have $c_d\sim Kd^{-\gamma} e^{-\lambda d}$
with some positive constants $K, \gamma, \lambda$; that is, the degree
distribution has a polynomial decay with exponential cut-off. On the
contrary, Chung,  Lu, Dewey, and Galas \cite[2003]{chung} claimed that for
$\beta=0$, when we do not have any extra edges, the asymptotic degree
distribution exists, and $(c_d)$ is decaying polynomially. None of
these papers contained a mathematically rigorous proof.  

Bebek, Berenbrink, Cooper, Friedetzky, Nadeau  and Sahinalp
\cite[2006]{bebek} disclaimed the above mentioned results of
\cite{pastor} and \cite{chung}. In the latter case, they showed that
the fraction of  isolated vertices (that have no edges) increases with
time in the pure duplication model, where $\beta=0$. They modified the
model to avoid singletons by adding a fixed number of edges to the new
vertex, chosen uniformly at random. They assumed without any proof
that the asymptotic degree distribution exists, and they claimed that
it is decaying polynomially.  

Hamdi, Krishnamurthy and Yin \cite[2013+]{hamdi} present a model where
the probabilities of adding a duplicated edge depends on the state of
a hidden Markov chain. Polynomial decay is stated for the limit of
the mean of degree distribution. We also mention the somewhat
different model of Jordan \cite[2011]{jordan}, and the duplication
model of Cohen,  Jordan and Voliotis \cite[2010]{cohen}, where the
duplicated vertex is chosen not uniformly, but with probabilities
proportional to the actual degrees.   

In our paper, we present a simple random graph model based on the
duplication of  a vertex chosen uniformly at random, and the erasure
of the edges of another vertex also chosen uniformly at random. We prove
that for all $d$, the proportion of vertices of degree $d$ tends to
some $c_d$ with probability $1$ as the number of steps goes to
infinity. Here $c_d$ is a positive number; we will formulate it as an
integral, and then we will determine the asymptotics of the sequence $(c_d)$ as
$d\rightarrow \infty$, showing that it has a stretched exponential
decay. Hence this model does not have the scale free property. We use
methods of martingale theory for proving almost sure convergence, and
generating function and Taylor series techniques for deriving the 
integral representation and the asymptotics of the sequence $(c_d)$.

\section{Definition of the model and main results}
Our model has two different versions. Both  of them start with a single
vertex. The graph evolves in discrete time steps; each step has a
duplication and an erasure part. At each step a new vertex will be
born; therefore the number of vertices after $n$ steps is $n+1$. The
graph is always a simple graph; it has neither multiple edges nor
loops. At each step we do the following.

\textit{Version 1}. We choose two (not necessarily different) old vertices
independently, uniformly at random. Then the new vertex is added to
the graph; we connect it to the first vertex and to all its
neighbors. After that we delete all edges emanating from the second
old vertex we have selected, with the possible exception that edges of
the new vertex cannot be deleted. 

\textit{Version 2}. We choose two (not necessarily different) old vertices
independently, uniformly at random. The new vertex is connected to the
first one and to its neighbors. Then we delete all edges of the
second vertex without any exceptions. 

That is, the new edges are protected in the erasure part of the same
step in version 1, but they might be deleted immediately in version 2.  
We will see that the version 2 graph has a simple structure that
enables us to describe its asymptotical degree distribution. Then,
using this and a coupling of the two models, we can prove similar
results for version 1. 

Let us remark that the presence of deletion makes the analysis more
difficult than in the usual recursive graph models, since it causes
intensive fluctuation in the model's behavior. 

Our model is a kind of coagulation--fragmentation one: the
effect of duplication is coagulation, and deletion results
fragmentation. Coagu\-lation--fragmentation models are frequently used
in several areas, see e.g. \cite{dong}. R\'ath and T\'oth applied these models 
for random graph models \cite{bal}, namely, for the Erd\H os--R\'enyi model, which is completely 
different from ours. 
  
The basic property of version 2 is that the evolving graph   
always consists of separated complete graphs. That is, it is a
disjoint union of cliques. Within a component, every pair of vertices
is connected, and there are no edges between the components. Indeed,
we start from a single vertex, which is a clique of size one, and both
duplication and  erasure make cliques from cliques. Moreover, it is
easy to see that if we start the model with an arbitrary graph, all
edges of the initial configuration are deleted after a while, and
after that the graph will consist of separated cliques. So the initial
configuration does not make any difference asymptotically.   

We may formulate the second version as follows. At each step we choose
two components independenty such that the probability that a given
clique is chosen is proportional to its size. The new vertex is
attached to the first clique, so its size is increased by $1$; the size
of the secondly chosen clique is decreased by $1$, and an isolated
vertex (the deleted one) comes into existence. Note that if we choose
an isolated vertex to be deleted, then it remains isolated.  

This structure of version 2 makes it easier to handle, as the number
of $d$-cliques does not vary so vehemently as the number of degree $d$
vertices; the fluctuation is bounded by $2$. This will lead 
to the description of the asymptotic degree distribution of version 1 in 
an almost sure sense. Our main results are the following.

\begin{theorem}\label{dupt1}
Denote by $X[n,d]$ the number of vertices of degree $d$ after $n$
steps in version 1. Then 
\[
\frac{X[n,d]}{n+1}\rightarrow c_d
\]
holds almost surely as $n\rightarrow\infty$, where $(c_d)$ is a
sequence of positive numbers satisfying 
\begin{equation} \label{duprekurzio}
c_0= \frac{1+c_1}{3}; \quad c_d=\frac{d+1}{2d+3}\bigl(c_{d-1}+c_{d+1}\bigr) \
\ \ \  (d\geq 2). 
\end{equation}
\end{theorem} 

For the asymptotic analysis we first present an integral representation for
the limiting sequence $(c_d)$. As a corollary, we get that the sum of this 
sequence is $1$; it is really a probability distribution.

\begin{theorem}\label{int} For the sequence $(c_d)$ of Theorem
  $\ref{dupt1}$ we have  
\[
c_d=(d+1)\int_0^\infty\frac{y^d e^{-y}}{(1+y)^{d+2}}\,dy \qquad (d\geq 0), 
\]
and $\sum_{d=0}^{\infty} c_d=1$.
\end{theorem}

Using this formula we can derive the asymptotics of $c_d$.

\begin{theorem}\label{asymp}For the sequence $(c_d)$ of Theorem
  $\ref{dupt1}$ we have  
\[
c_d\sim (e\pi)^{1/2}\,d^{1/4}\,e^{-2\sqrt{d}},\quad as\ d\to\infty.
\]
\end{theorem}

Our model is invented to ensure high degree clustering. Finally, let
us quantify this property. 

The local clustering coefficient of a vertex of degree $d$ is defined
to be the fraction of connections that exist between the $\binom{d}{2}$ 
pairs of neighbors (meant $0$ when $d<2$). Watts and Strogatz
\cite{WS} define the clustering coefficient of the whole graph as the
average of the local clustering coefficients of all the vertices. Let
us call this quantity the average clustering coefficient. Another
possibility for a such a measure is the ratio of $3$ times the number
of triangles divided by the number of connected triplets (paths of
length $2$), see \cite{hofstad}. This version is sometimes called
transitivity; we will refer to it as the global clustering coefficient. 

Since the graph in version 2 is consists of disjoint cliques, its global 
clustering coefficient is obviously $1$, while the
average clustering coefficient is equal to the proportion of vertices
with degree at least $2$. By Theorem \ref{dupt1} it converges to
$1-c_0-c_1=2-4c_0$ almost surely, as $n\to\infty$. We note that the
limit is equal to $0.38538\dots$ by Theorem \ref{int}. These
results can be transferred to version 1.

\begin{theorem}\label{clust}
In version 1, the global clustering coefficient converges to $1$,
and the average clustering coefficient to $1-c_0-c_1$, almost
surely, as $n\to\infty$. 
\end{theorem}

The high clustering property of our model shows that is is a so-called small-world graph \cite{WS}.

\section{Proofs}

\subsection*{Preliminaries.} First we formulate the lemma from
martingale theory that we will use several times and whose proof can
be found in \cite{publ}.

\begin{lemma}\label{publ1}
Let $(\mathcal F_n)$ be a filtration, $(\xi_n)$ a nonnegative adapted
process. Suppose that 
\begin{equation}\label{lemmafelt}
\mathbb \mathbb{E}\bigl((\xi_n-\xi_{n-1})^2\bigm|\mathcal
F_{n-1}\bigr)=O\left( n^{1-\delta}\right)
\end{equation}
holds with some $\delta>0$. Let $(u_n)$, $(v_n)$ be nonnegative
predictable processes such that $u_n<n$ for all $n\geq 1$. Finally,
let $(w_n)$ be a regularly varying sequence of positive numbers with
exponent $\mu\ge 0$.

$(a)$ Suppose that 
\[
\mathbb \mathbb{E}(\xi_n\mid\mathcal F_{n-1})\le
\Bigl(1-\dfrac{u_n}{n}\Bigr)\xi_{n-1}+v_n,
\]
and $\lim_{n\rightarrow \infty} u_n=u$, $\limsup_{n\rightarrow \infty}
v_n/w_n\le v$ with some random variables $u>0,\ v\geq 0$. Then  
\begin{equation*}
\limsup_{n\rightarrow\infty}\frac{\xi_n}{nw_n}\le \frac{v}{u+\mu+1}
\quad a.s.
\end{equation*} 

$(b)$ Suppose that 
\[
\mathbb \mathbb{E}(\xi_n\mid\mathcal F_{n-1})\ge
\Bigl(1-\dfrac{u_n}{n}\Bigr)\xi_{n-1}+v_n,
\]
and $\lim_{n\rightarrow \infty}u_n=u$, $\liminf_{n\rightarrow \infty}
v_n/w_n\ge v$ with some random variables $u>0,\ v\geq 0$. Then  
\begin{equation*}
\liminf_{n\rightarrow\infty}\frac{\xi_n}{nw_n}\ge \frac{v}{u+\mu+1}
\quad a.s.
\end{equation*}
\end{lemma}

\subsection*{Asymptotic degree distribution in version 2.} 
Recall that in this case the graph is always disjoint union of
complete graphs. 

First we prove the following analogue of Theorem \ref{dupt1}.

\begin{prop}\label{dupall1}
Denote by $Y[n,k]$ the number of cliques of size $k$ after $n$ steps
in version 2. Then for all positive integers $k$ we have 
\[
\frac{Y[n,k]}{n}\rightarrow y_k \quad {\text\ almost\ surely\  as\ }
n\rightarrow \infty,
\]
where $(y_k)$ is a sequence of positive numbers satisfying
\begin{equation}\label{dupe0}
y_1=\frac{1+2y_2}{3}, \qquad
y_k=\frac{(k-1)y_{k-1}+(k+1)y_{k+1}}{2k+1} \quad (k\geq 2).
\end{equation}
\end{prop}

Note that \eqref{dupe0} (as well as equation \eqref{duprekurzio}) is
not a recursion. This prevents us proceeding simply in the usual,
direct way, with induction over $k$.

\textbf{Proof.} For $n=0$ we have $Y[0,1]=1$, all the other ones are equal
to zero. The total number of vertices is $n$ after $n-1$
steps. Let $\mathcal F_n$ denote the $\sigma$-field generated by the
first $n$ steps.

We enumerate the events that can happen to the cliques of different
sizes at a step. 
 
At the $n$th step an isolated vertex may become 
\begin{itemize}
\item a clique of size 2 (increased but not decreased) with
  probability $\frac1n \big(1-\frac 1n\big)$; 
\item an isolated vertex (any other cases). 
\end{itemize}

A clique of size $k\geq 2$ may become a clique of size 
\begin{itemize}
\item  $k-1$ (not increased but decreased) with probability $\frac kn
\big(1-\frac kn\big)$;
\item $k+1$ (increased but not decreased) with probability $\frac kn
  \big(1-\frac kn\big)$;
\item $k$ (any other cases). 
\end{itemize}

The deleted vertex will be a new isolated one unless one of them is
chosen for erasure  but not for duplication, which has probability
$\frac1n \big(1-\frac 1n\big)$ for each of them. 

Putting this together with the fact that the random choices are
independent and probabilities are proportional to clique sizes, we can
compute the conditional expectation of $Y[n,k]$ with respect to
$\mathcal F_{n-1}$, which is the $\sigma$-field generated by the
first $n-1$ steps.

\begin{align*}
E( Y[n,1]\vert \mathcal
F_{n-1})&=Y[n-1,1]\bigg[1-\frac{1}{n}\bigg(1-\frac{1}{n}\bigg
)-\frac{1}{n}\bigg(1-\frac{1}{n}\bigg )\bigg ]\\
&\quad +1+ Y[n-1,2]\cdot \frac{2}{n}\bigg(1-\frac 2n\bigg); \\
E(Y[n,k]\vert \mathcal F_{n-1})&=Y[n-1,k]\bigg [1-2\cdot \frac
kn\bigg(1-\frac{k}{n}\bigg)\bigg]\\
&\quad +Y[n-1,k-1] \cdot \frac{k-1}{n} \bigg (1-\frac{k-1}{n}\bigg)\\
&\quad+Y[n-1,k+1]\cdot \frac{k+1}{n}\bigg(1-\frac{k+1}{n}\bigg) 
\quad (k\geq 2).
\end{align*}

Let $A_k=\liminf_{n\rightarrow\infty} \frac{Y[n,k]}{n}$ and
$B_k=\limsup_{n\rightarrow\infty} \frac{Y[n,k]}{n}$ 
for  $k\geq 1$. It is clear that $0\leq A_k\leq B_k\leq 1$ holds for
these random variables. 

We will give a sequence of lower bounds for $(A_k)$, and similarly, a
sequence of upper bounds for $(B_k)$; the we will show that
their limits are equal to each other. First, let $a_k^{(0)}=0$ for
$k\geq 1$. Having constructed the sequence $(a_k^{(j)})_{k\geq 1}$, we
define 
\begin{equation}\label{dupe1}
a_1^{(j+1)}=\frac{1+2a_2^{(j)}}{3}\,, \quad
a_k^{(j+1)}=\frac{(k-1)a_{k-1}^{(j)}+(k+1)a_{k+1}^{(j)}}{2k+1} \quad
(k\geq 2).
\end{equation}
We get $a_k^{(j)}$ recursively for every $k\geq 1$ and $j\geq 1$. 

We prove by induction on $j$ that $a_k^{(j)}\leq A_k \ (k\geq
1)$. Since  $Y[n,k]\geq 0$, this is clear for $j=0$. Suppose that this
is satisfied for some $j$ for every $k$. For $k=1$ we apply Lemma
\ref{publ1} with 
\[
\xi_n={Y[n,1]}, \ u_n=2-\frac{2}{n}\rightarrow 2, \
v_n=1+Y[n-1,2]\cdot\frac{2}{n}\bigg(1-\frac 2n\bigg).
\]

Now $(\xi_n)$ is nonnegative adapted. $(u_n)$ and $(v_n)$ are clearly
nonnegative predictable sequences; we can choose $w_n=1$, $\mu=0$,
$u=2>0$ and finally, $v=1+2a_2^{(j)}\geq 0$ due to the induction
hypothesis. Note that at each step at most one of the isolated points 
vanishes and at most two may appear. Thus \eqref{lemmafelt} is clearly
satisfied. Lemma \ref{publ1} implies that 
\[
A_1=\liminf_{n\rightarrow \infty}
\frac{Y[n,1]}{n}=\liminf_{n\rightarrow \infty} \frac{\xi_n}{n} \geq
\frac{v}{u+1}=\frac{1+2a_2^{(j)}}{3}=a_1^{(j+1)}
\]
almost surely. 

Similarly, for $k\geq 2$, if we have $A_k\geq a_k^{(j)}$ for some
$j\geq 1$, we can choose 
\begin{align*}
\xi_n&={Y[n,k]}, \ u_n=2k-\frac{2k^2}{n}\rightarrow 2k,\\ 
v_n&=Y[n-1,k-1]\cdot\frac{k-1}{n}\bigg (1-\frac{k-1}{n}\bigg)\\
&\quad +Y[n-1,k+1]\cdot \frac{k+1}{n}\bigg(1-\frac{k+1}{n}\bigg),\\
v&=(k-1)a_{k-1}^{(j)}+(k+1)a_{k+1}^{(j)}.
\end{align*}
At each step at most three cliques are changed, which implies that
\eqref{lemmafelt} holds. Thus in this case from Lemma \ref{publ1} we
obtain that 
\[
A_k=\liminf_{n\rightarrow \infty} \frac{Y[n,k]}{n} \geq \frac{v}{u+1}=
\frac{(k-1)a_{k-1}^{(j)}+(k+1)a_{k+1}^{(j)}}{2k+1}=a_k^{(j+1)}
\]
almost surely.

By induction on $j$ we get that $A_k\geq a_k^{(j)}$ holds almost
surely for  $k\geq 1$ and $j\geq 0$.

Now we verify that for fixed $k$ the sequence $(a_k^{(j)})$ is
monotone increasing in $j$. Since $a_k^{(0)}=0$ for every $k$, from
equations \eqref{dupe1} it is clear that $a_k^{(1)}\geq
a_k^{(0)}$. Suppose that for some $j\geq 1$ we have $a_k^{(j)}\geq
a_k^{(j-1)}$ for every $k$. Then 
\begin{align*}
a_1^{(j+1)}&=\frac{1+2a_2^{(j)}}{3}\geq \frac{1+2a_2^{(j-1)}}{3}=
a_1^{(j)};\\
a_k^{(j+1)}&=\frac{(k-1)a_{k-1}^{(j)}+(k+1)a_{k+1}^{(j)}}{2k+1}\\
&\geq \frac{(k-1)a_{k-1}^{(j-1)}+(k+1)a_{k+1}^{(j-1)}}{2k+1}
=a_k^{(j)}
\end{align*}
follows from equations \eqref{dupe1}. Thus by induction on $j$ we get
that $a_k^{(j)}\geq a_k^{(j-1)}$ for  $k, j\geq 1$.

It is clear that the sequence $(a_k^{(j)})_{j\geq 0}$ is uniformly
bounded from above by $1$. Using monotonicity we can define 
\[
a_k= \lim_{j\rightarrow \infty }a_k^{(j)} \qquad (k\geq 1) .\]
From equation \eqref{dupe1} it follows that $(a_k)$ satisfies
\eqref{dupe0}, that is,  
\begin{equation*}
a_1=\frac{1+2a_2}{3}, \qquad
a_k=\frac{(k-1)a_{k-1}+(k+1)a_{k+1}}{2k+1} \quad (k\geq 2).
\end{equation*}
On the other hand, since $A_k\geq a_k^{(j)}$ for  $k\geq 1$ and $j\geq
0$, we have $A_k\geq a_k$ almost surely.

Similarly, we define $b_k^{(0)}=1$ for every $k$, and then
\[
b_1^{(j+1)}=\frac{1+2b_2^{(j)}}{3}, \qquad b_k^{(j+1)}=
\frac{(k-1)b_{k-1}^{(j)}+(k+1)b_{k+1}^{(j)}}{2k+1} \quad (k\geq 2).
\]
Using part (a) of Lemma \ref{publ1} it follows by induction on $j$
that $B_k\leq b_k^{(j)}$ holds almost surely.

In this case, for fixed $k$ the sequence $b_k^{(j)}$ is decreasing,
and for the limits $b_k=\lim_{j\rightarrow \infty} 
b_k^{(j)}$ we also have 
\begin{equation*}
b_1=\frac{1+2b_2}{3}, \qquad
b_k=\frac{(k-1)b_{k-1}+(k+1)b_{k+1}}{2k+1} \quad (k\geq 2).
\end{equation*}
In addition, $B_k\leq b_k$ almost surely.

By definition, $0\leq A_k\leq B_k\leq 1$ and $0\leq a_k\leq b_k\leq 1$
hold. Let $d_k=b_k-a_k\geq 0$ for all $k$. We have the same equations
for $(a_k)$ and $(b_k)$. This yields 
\begin{equation*}
d_1=\frac{2d_2}{3}, \qquad d_k=\frac{(k-1)d_{k-1}+(k+1)d_{k+1}}{2k+1} 
\quad (k\geq 2).
\end{equation*}
By rearranging we get that
\begin{equation}\label{dupe4}
d_2=\frac{3}{2}d_1, \qquad
d_{k+1}=\frac{(2k+1)d_{k}-(k-1)d_{k-1}}{k+1} \quad (k\geq 2).
\end{equation}

Suppose that $d_k\geq \frac{k+1}{k} d_{k-1}$ holds for some $k\geq
2$. (For $k=2$ this is true with equality.) Since $d_{k-1}$ is
nonnegative, $d_k\geq d_{k-1}$ also follows from this assumption. We
obtain from equation \eqref{dupe4} that 
\[
d_{k+1}\geq \frac{(k+2)d_k}{k+1}.
\] 
Therefore this inequality holds for every $k$. 

This implies that $d_k\geq (k+1) d_1$ for every $k$. Since $0\leq
d_k=b_k-a_k\leq 1$, it follows that $d_1=0$.

From \eqref{dupe4} we obtain that $d_k=0$ for
all $k$, which implies that $a_k=b_k$. Since these were the lower and
upper bounds for the limit inferior and limit superior of
$\frac{Y[n,k]}{n}$, we get that the latter must converge almost surely
as $n\to\infty$, and the limits satisfy \eqref{dupe0}. \qed 

\begin{corollary}\label{dupk2}
In version 2, the proportion of vertices of degree $d$ tends to
$c_d$ satisfying \eqref{duprekurzio} almost surely as $n\rightarrow
\infty$.
\end{corollary}

\textbf{Proof.} For a fixed $d$ we have $d+1$ vertices of degree $d$
in each clique of size $k=d+1$. Therefore for the proportion of
vertices of degree $d$ tends to $(d+1)y_{d+1}$ by Proposition
\ref{dupall1}. From equations \eqref{dupe0} we obtain that 
\[
c_0=y_1=\frac{1+2y_2}{3}=\frac{1+c_1}{3};
\]
\[
c_d=(d+1)y_{d+1}=\frac{d+1}{2d+3}\bigl(c_{d-1}+c_{d+1}\bigr)
\quad (d\geq 2).
\] 
\qed

\subsection*{Asymptotic degree distribution in version 1.}

When proving the results for version 2 we essentially used the
property that the graphs consists of disjoint union of cliques: at
most three of the cliques may change at a step, but the number of
vertices whose degree is changed is not bounded uniformly. However,
we can push through the results by a kind of coupling of versions 1
and 2.  

{\bf Proof of Theorem \ref{dupt1}.} 
Both in versions 1 and 2 two old vertices are selected with
replacement, independently, uniformly at random. Thus we can couple
the models such that the selected vertices are the same in all
steps. The duplication part is the same in the two versions. The
difference is in the deletion: in version 1, the edges of the new vertex
cannot be deleted. So in version 1, we do the following. In the
deletion part, we colour an edge red if it is saved in version
2. That is, if it connects the new vertex with the old vertex to be
deleted. In the duplication part, copies of red edges are also red: if
there is a red edge between the duplicated vertex and one of its
neighbors, then the new edge connecting this neighbor to the new
vertex is also red. All other new edges are originally black, but they
may turn red in the deletion part of the same step.  

The colouring is defined in such a way that the graph sequence of the
black edges is a 
realization of version 2. Indeed, edges turning red are deleted and
hence the copies of them does not appear in this model, but all other
edges are black.   

Our goal is to prove that the number of vertices having red edges
divided by $n$ tends to zero almost surely. This implies that the
results of Corollary \ref{dupk2} holds for version 1 as well.  

First we need an upper bound for the total number of edges. 

\begin{lemma}\label{lemma2}
Denote by $S_n$ the number of edges (both black and red ones) after
$n$ steps in version 1. Then for all $\varepsilon>0$ we  
have $S_n=O\bigl(n\log^{1+\varepsilon}n\bigr)$ with probability $1$.
\end{lemma}
{\bf Proof.} Let
$\delta_n=S_n-S_{n-1}$. As before, $\mathcal F_n$ denotes the
$\sigma$-field generated by the first $n$ steps, and  $X[n,d]$ is
the number of vertices of degree $d$ after $n$ steps. Let $U_n$,
resp. $V_n$, denote the degree of the old vertex selected for
duplication, resp. deletion, at step $n$. The new vertex is connected
to the duplicated one with an edge that cannot be deleted; this
increases the number of edges by $1$ for sure. Thus,
$\delta_n=U_n-V_n+1$. Clearly, $U_n$ and $V_n$ are conditionally
i.i.d. with respect to $\mathcal F_{n-1}$, hence
$S_n-n=\sum_{j=1}^n(\delta_j-1)$ is a zero mean
martingale. Consequently, $ES_n=n$ for every $n$.  

Clearly, 
\[
E\bigl(|\delta_n-1|\bigm|\mathcal F_{n-1}\bigr)\le 2 
E(U_n\vert\mathcal F_{n-1})=\sum_{d=0}^n \frac{X[n-1,d]}{n}\,d=
\frac{2S_{n-1}}{n}\,.
\]
Hence
\[
E\left(\sum_{n=2}^\infty\frac{|\delta_n-1|}{n\log^{1+\varepsilon}n}
\right)<\infty,
\]
therefore the series
\[
\sum_{n=2}^\infty\frac{\delta_n-1}{n\log^{1+\varepsilon}n}
\]
is convergent with probability $1$. Then Kronecker's lemma \cite[Lemma
IV.3.2]{sh} implies that
\[
\frac{S_n-n}{n\log^{1+\varepsilon}n}\to 0\quad\text{a.s.}
\]
as $n\to\infty$.\qed

Now we will colour some of the vertices red in such a way that the
remaining black vertices cannot have any red edges. We will be able
to give an upper bound for the number of red vertices.  

At the duplication step the new vertex becomes red if and only if the
duplicated vertex is red. If this old vertex is black and has no red
edges, the same holds for the new vertex at the moment. After that, if
there is an edge between the new vertex and the deleted one, this edge
may turn red, as we defined before. We colour both endpoints of this
new red edge red. On the other hand, if the old vertex chosen for
deletetion loses all its edges, then its new colour will be black. Note
that black vertices still have only black edges, but it may happen
that an old vertex has only one red edge which is deleted, because
its other endpoint is chosen for deletion; in this case the vertex
stays red without having any red edges.

The proof continues with giving an upper bound for the number of 
red vertices.

\begin{lemma}\label{lemma3}
Denote by $Z_n$ the number of red vertices after $n$ steps. Then for all 
$\varepsilon>0$ we have $Z_n=O(\log^{2+\varepsilon}n)$ almost surely.
\end{lemma}
{\bf Proof.}
At each step, every old
vertex has the same probability to be duplicated or deleted. If a red
vertex is duplicated, then the new vertex becomes red; if it is
deleted, then $Z_n$  decreases by 1 unless the deleted vertex is
connected to the new one which turns this edge red. Therefore without
the exceptional new red edge, the conditional expectation of $Z_n$
with respect to $\mathcal F_{n-1}$ would be equal to $Z_{n-1}$. The
deleted vertex and the new one are connected if and only if the
deleted and duplicated vertices are the same or they are connected to
each other. Since we did sampling with replacement, the probability of
the first event is $1/n$; while the probability of the second event is
$2S_{n-1}/n^2$. In the first case, the new vertex is red originally, but
the other one stays red instead of turning back to black when deleted;
$Z_n$ is increased by an extra 1. In the other case, both endpoints of
the edge turning red may be  red vertices in addition. To sum up, we
obtain that
\[
E(Z_n\vert \mathcal F_{n-1})\leq Z_{n-1}+\frac 1n+
4\cdot\frac{S_{n-1}}{n^2}.
\]  

We set $\eta_n=Z_n-Z_{n-1}$. With this notation 
\begin{equation}\label{Ed_n}
E(\eta_n\vert \mathcal F_{n-1})\leq \frac 1n+ 4\cdot\frac{S_{n-1}}{n^2}
\end{equation}

We have already shown that $ES_{n-1}=n-1$, hence $E\eta_n\le 5/n$, and
$EZ_n=O(\log n)$. 

Note that the number of red vertices cannot change
by more than three at a single step, because if an old vertex is neither
deleted, nor duplicated, it cannot be coloured red. Hence $|\eta_n|\leq
3$ for all $n$. 
Moreover, we can give an upper bound on the probability 
that the number of red vertices changes at step $n$. Namely, 
it can change only if 
\begin{itemize}
\item we duplicate and delete the same vertex; this has (conditional)
  probability $1/n$.  
\item the duplicated and the deleted vertices are connected to each
  other; this has probability $2S_{n-1}/n^2$, because   
 there are $S_{n-1}$ edges.
 \item a red vertex is duplicated; this has probability $Z_{n-1}/n$. 
 \item a red vertex is deleted; this has probability $Z_{n-1}/n$.
\end{itemize}

Thus 
\begin{equation}\label{Vd_n}
P(Z_n\neq Z_{n-1}\vert \mathcal F_{n-1})\leq 
\frac 1n +2\frac{S_{n-1}}{n^2}+2\frac{Z_{n-1}}{n},
\end{equation}
therefore
\[
E|\eta_n|\le 3P(Z_n\ne Z_{n-1})=O\Bigl(\frac{\log n}{n}\Bigr),
\]
which implies that
\[
E\left(\sum_{n=2}^\infty\frac{|\eta_n|}{\log^{2+\varepsilon}n}\right)
<\infty.
\]
The proof can be completed by the help of Kronecker's lemma, just like in
the proof of Lemma \ref{lemma2}.\qed

Now we can finish the proof of Theorem \ref{dupt1}.

The total number of vertices is $n+1$ after $n$ steps, hence the
proportion of red vertices converges to 0 almost surely as
$n\rightarrow \infty$.  Since we defined the colours in such a way
that red edges are exactly the edges that are present in version 1 but
are not present in version 2, and only red vertices may have red
edges, it follows that the proportion of vertices having  different
degree in the two versions converges to 0. Corollary \ref{dupk2}
states that for every $d$ the proportion of vertices of degree $d$ in
version 2 converges almost surely to $c_d$. Now the same follows for
version 1, which is the statement of Theorem \ref{dupt1}.\qed

\begin{remark}
We could have given an upper bound for the conditional expectation of
the number of red edges. The advantage of using red vertices is the
uniform bound on the total change in their number; there is no such
bound for the change in the number of red edges. 
\end{remark}

\begin{remark}
It follows that version 1 has a quite specific structure: it consists
of cliques that are connected with relatively few edges (those are
coloured red). An edge can be red only if both its endpoints are red,
hence Lemma \ref{lemma3} gives an $O\bigl(\log^{4+\varepsilon}n\bigr)$ 
bound for the number of red edges.
\end{remark}
This is not sharp; however, the estimates of Lemmas \ref{lemma2} and
\ref{lemma3} can be further improved, which might be, as pointed out
above, of independent interest. Thus, before turning to the proof of
Theorem \ref{int}, we present the following improvement.
\begin{prop}\label{sharp}
$S_n\sim n$, and $Z_n=O\bigl(\log^{1+\varepsilon}n\bigr)$ for every
$\varepsilon >0$ almost surely, as $n\to\infty$.
\end{prop}
\textbf{Proof.}
First we give a crude bound for the maximal degree $M_n=\max\{d:
X[n,d]>0\}$.  According to Lemma \ref{lemma2}, $S_n=O\bigl(n\log^{1+\varepsilon}n\bigr)$ also holds for the 
number of edges in version 2.
Since a clique of size $k$ contains $\binom{k}{2}$ edges,  it
follows that the size of the maximal clique is
$O\bigl(n^{1/2+\varepsilon}\bigr)$. The same holds for the maximal
degree in version 2; and, by Lemma \ref{lemma3}, in version 1, too.
Thus $M_n=O\bigl(n^{1/2+\varepsilon}\bigr)$ for every
$\varepsilon>0$. 

Next, consider the martingale $S_n-n=\sum_{j=1}^n(\delta_j-1)$ from
the proof of Lemma \ref{lemma2}. In order to prove that
$S_n-n=o(\gamma_n)$ for a positive increasing predictable sequence
$(\gamma_n)$ it is sufficient to show that 
\[
\sum_{n=1}^\infty\gamma_n^{-2}\,E\bigl((\delta-1)^2\bigm|\mathcal F_{n-1}
\bigr)<\infty
\]
with probability $1$ \cite[Theorem VII.5.4]{sh}. To this end we need
to estimate the conditional variance of the martingale differences.
\begin{multline*}
\mathrm{Var}(\delta_n-1\vert\mathcal F_{n-1})=
2\mathrm{Var}(U_n\vert\mathcal F_{n-1})\le
2E(U_n^2\vert\mathcal F_{n-1})\\
=2\sum_{d=1}^n \frac{X[n-1,d]}{n}\,d^2\le \frac 2n\,M_{n-1}
\sum_{d=1}^n X[n-1,d]\,d\\
=\frac 2n\,M_{n-1}S_{n-1}=
O\bigl(n^{1/2+\varepsilon}\bigr),
\end{multline*}
for every positive $\varepsilon$. Hence
\[
\sum_{n=1}^\infty \frac{E\bigl((\delta-1)^2\bigm|\mathcal F_{n-1}
\bigr)}{n^{3/2+\varepsilon}}<\infty,
\]
implying 
\[
S_n-n=o\bigl(n^{3/4+\varepsilon}\bigr)
\]
Thus $S_n\sim n$ a.s., indeed.

Finally, let us consider the martingale
$\zeta_n=\sum_{j=1}^n\bigl(\eta_j-E(\eta_j|\mathcal F_{j-1})\bigr)$,
 where $\eta_n=Z_n-Z_{n-1}$, and 
derive an upper bound for the conditional variance of the differences.
Keeping in mind that $|\eta_n|\le 3$ and using \eqref{Vd_n} we have
\begin{multline*}
E\bigl((\zeta_n-\zeta_{n-1})^2\bigm| \mathcal F_{n-1}\bigr)=
\textrm{Var}(\eta_n\vert \mathcal F_{n-1})\\
\leq E((Z_n-Z_{n-1})^2\vert \mathcal F_{n-1})
\le 9 P(Z_n\neq Z_{n-1}\vert \mathcal F_{n-1})\\
\le 9\bigg ( \frac 1n +2\frac{S_{n-1}}{n^2}+2\frac{Z_{n-1}}{n}\bigg)
=O\Bigl(\frac{1+Z_{n-1}}{n}\Bigr).
\end{multline*}

Now suppose that  $Z_n=O(\log^\alpha n)$ is satisfied for some
$\alpha>0$. Then 
\[
E\bigl((\zeta_n-\zeta_{n-1})^2\bigm| \mathcal F_{n-1}\bigr)=
O\Bigl(\frac{\log^{\alpha}n}{n}\Bigr),
\]
hence
\[
\sum_{n=2}^\infty \frac{E\bigl((\zeta_n-\zeta_{n-1})^2\bigm| 
\mathcal F_{n-1}\bigr)}{\log^{\alpha+1+\varepsilon}n}<\infty
\]
with probability $1$. Again, by \cite[Theorem VII.5.4]{sh} we have
\begin{equation}\label{zeta}
\zeta_n=o\bigl(\log^{(\alpha+1)/2+\varepsilon}\bigr)\quad\text{a.s.}
\end{equation}
for every positive $\varepsilon$.

Clearly,
\[
Z_n=\sum_{j=1}^n \eta_j=\zeta_n+\sum_{j=1}^n E(\eta_j|\mathcal F_{j-1}),
\]
where the last sum can be estimated by the help of \eqref{Ed_n} in the
following way. Since $S_{n-1}\sim n$, we have $E(\eta_n|\mathcal
F_{n-1})=O(1/n)$, hence 
\[
\sum_{j=1}^n E(\eta_j|\mathcal F_{j-1})=O(\log n).
\]
This, combined with \eqref{zeta} gives that $Z_n=
O\bigl(\log^{(\alpha+1)/2+\varepsilon}\bigr)$ holds almost
surely for all $\varepsilon>0$. By Lemma \ref{lemma3} we can start
from $\alpha=2+\varepsilon$, and repeating the argument we finally end
up with the a.s. estimation $Z_n=O\bigl(\log^{1+\varepsilon}\bigr)$,
for all $\varepsilon>0$. \qed 
\bigskip

{\bf Proof of Theorem \ref{int}.}
Let $G(z)$ denote the generating function of the sequence $(c_d)$,
that is,
\[
G(z)=\sum_{d=0}^{\infty}c_dz^d,\quad |z|\le 1.
\]
Multiplying equation $(d+1)(c_{d-1}+c_{d+1})=(2d+3)c_d$ by $z^d$, then
summing up from $d=1$ to $\infty$ and using that $c_0=(1+c_1)/3$, we
obtain an inhomogeneous linear differential equation for $G(z)$.
\[
(1-z)^2G^{\prime}(z)=(3-2z)G(z)-1,\quad G(0)=c_0.
\]
Solving this equation we get the following expression
\[
G(z)=\frac{c(z)}{(1-z)^2}\,\exp\Bigl(\frac{z}{1-z}\Bigr),
\]
where
\[
c(z)=c_0-\int_0^z\exp\Bigl(-\frac{y}{1-y}\Bigr)\,dy.
\]
Since $G(1)=\sum_{d=0}^{\infty} c_d\le 1$, it follows that
\[
c_0=\int_0^1\exp\Bigl(-\frac{y}{1-y}\Bigr)\,dy,
\]
hence, via the substitution $x=1-y$,
\[
c(z)=\int_z^1\exp\Bigl(-\frac{y}{1-y}\Bigr)\,dy=
\int_0^{1-z}\exp\Bigl(1-\frac{1}{x}\Bigr)\,dx.
\]
Thus we have
\[
G(z)=\int_0^{1-z}\exp\Bigl(1-\frac{1}{x}\Bigr)\,dx\ \frac{1}{(1-z)^2}\,
\exp\Bigl(\frac{z}{1-z}\Bigr),
\]
from which, by substituting $y=\frac 1x-\frac{1}{1-z}$, we obtain
\begin{multline}\label{genfv}
G(z)=\int_0^\infty\frac{e^{-y}}{(1+(1-z)y)^2}\ dy\\
=\int_0^\infty\frac{e^{-y}}{(1+y)^2\bigl(1-z\,
\frac{y}{1+y}\bigr)^2}\ dy\\
=\int_0^\infty\sum_{d=0}^\infty(d+1)\,\frac{z^dy^d\,e^{-y}}
{(1+y)^{d+2}}\ dy\\
=\sum_{d=0}^\infty z^d\,(d+1)\int_0^\infty\frac{y^d e^{-y}}
{(1+y)^{d+2}}\ dy,
\end{multline}
completing the proof of the first statement of the theorem. 

In addition, note that the first equality of \eqref{genfv} immediately 
implies that $\sum_{d=0}^{\infty} c_d=G(1)=1$.
\qed

{\bf Proof of Theorem \ref{asymp}.}
In order to approximate the integral of Theorem \ref{int} we first
analyse the behavior of the integrand around the point where it 
attains its maximum. Let
\[
y_d=\arg\max\frac{y^d e^{-y}}{(1+y)^{d+2}}=\arg\max f(y),
\]
where
\[
f(y)=d\log y-(d+2)\log(1+y)-y.
\]
Clearly,
\begin{align*}
f^\prime(y)&=\frac dy-\frac{d+2}{y+1}-1=-\frac{y^2+3y-d}{y(y+1)}\,,\\
f''(y)&=-\frac{d}{y^2}+\frac{d+2}{(y+1)^2}=\frac{2y^2-2dy-d}{y^2(y+1)^2}
\,,\\
f'''(y)&=\frac{2d}{y^3}-\frac{2(d+2)}{(y+1)^3}\,.
\end{align*}

Since $y_d$ satisfies $f^\prime(y_d)=0$, we get that
\[
y_d=-\frac 32 + \sqrt{d+\frac 94}=\sqrt{d}-\frac 32+o(1).
\]
Let us write $y$ in the form $y=y_d+y_d^{1/2}t$. Then
\[
g(t):=f(y)-f(y_d)=\frac{y_d}{2}\,f''(y_d+\theta y_d^{1/2}t)\,t^2,
\]
where $\theta=\theta(d,t)$ belongs to the interval $[0;1]$.
For every fixed $t$
\[
f''(y_d+\theta y_d^{1/2}t)\sim -2y_d^{-1},
\]
thus $g(t)\to -t^2$ as $d\to\infty$. 
Moreover, for $y\le y_d$, that is, for $y_d^{1/2}\le t\le 0$ we have
$f'(y) \ge 0$. Thus  $d/y - (d+2)/(y+1) > 0$ holds, and after
rearranging we get that $(d+2)/d < (y+1)/y$. This yields that $(d+2)/d
< (y+1)^3/y^3$ is satisfied, which implies that $f'''(y)\ge 0$. Hence 
\[
g(t)\le\frac{y_d}{2}\,f''(y_d)\,t^2=a_d\,t^2,
\]
where $a_d\to -1$, as $d\to\infty$. On the other hand, let $y_d\le
y\le \frac 32\,y_d$, that is, $0\le t\le \frac 12\,y_d^{1/2}$. In this
domain $f'''$ is increasing, hence $f'''(y)\le f'''(y_d)\sim
6dy_d^{-4}\sim 6d^{-1}$. Thus, 
\begin{multline*}
g(t)\le\frac{y_d}{2}\,f''(y_d)\,t^2+\frac
16\,y_d^{3/2}f'''(y_d)\,t^3\\
\le\left(\frac{y_d}{2}\,f''(y_d)+\frac{y_d^2}{12}\,f'''(y_d)\right)t^2
=b_dt^2,
\end{multline*}
where $b_d\to -1/2$, as $d\to\infty$.

Thus, by the dominated convergence
theorem, 
\begin{multline*}
\int_0^{3y_d/2}e^{f(y)}\,dy=y_d^{1/2}\int_{-y_d^{1/2}}^{\frac 12\,y_d^{1/2}} 
\exp\bigl(f(y_d)+g(t)\bigr)\,dt\\
\sim y_d^{1/2}\exp\bigl(f(y_d)\bigr)\int_{-\infty}^{+\infty}
\exp\bigl(-t^2\bigr)\,dt=\sqrt{\pi}\,y_d^{1/2}
\exp\bigl(f(y_d)\bigr).
\end{multline*}

Here 
\[
f(y_d)=-2\log y_d-(d+2)\log\Bigl(1+\frac{1}{y_d}\Bigr)-y_d,
\]
and
\begin{align*}
(d+2)\log\Bigl(1+\frac{1}{y_d}\Bigr)&=(d+2)\Bigl(\frac{1}{y_d}-
\frac{1}{2y_d^2}\Bigr)+o(1)\\
&=y_d+\frac{(d+2)(2y_d-1)-2y_d^3}{2y_d^2}+o(1)\\
&=y_d+\frac{(y_d^2+3y_d+2)(2y_d-1)-2y_d^3}{2y_d^2}+o(1)\\
&=y_d+\frac{5y_d^2+y_d-2}{2y_d^2}+o(1),
\end{align*}
where we used that $y_d^2+3y_d=d$. Thus,
\[
f(y_d)=-2\log y_d-2y_d-\frac 52+o(1)=-2\log y_d-2\sqrt{d}+\frac
12+o(1).
\]

Finally, 
\begin{align*}
\int_{3y_d/2}^\infty e^{f(y)}\,dy&\le \bigl(2y_d\bigr)^{-2}
\int_{3y_d/2}^\infty \Bigl(1-\frac{1}{1+y}\Bigr)^{\!d}\,e^{-y}\,dy\\
&\leq\bigl(2y_d\bigr)^{-2}\int_{3y_d/2}^\infty \exp\Bigl(-\frac{d}{y+1}
-y\Bigr)\,dy.
\end{align*}
The exponent on the right-hand side can be estimated with the help of
the AM--GM inequality as follows.
\[
-\frac{d}{y+1}-y=-\frac{d}{y+1}-\frac{y+1}{2}-\frac{y-1}{2}\le
-\sqrt{2d}-\frac{y-1}{2},
\]
hence
\[
\int_{3y_d/2}^\infty e^{f(y)}\,dy\le \bigl(2y_d\bigr)^{-2}\,
\exp\Bigl(-\sqrt{2d}+\frac 12-\frac 34\,y_d\Bigr)=
o\Bigl(y_d^{-2}\,\exp\bigl(-2\sqrt{d}\bigr)\Bigr).
\]
From all these we obtain that
\[
c_d=(d+1)\int_0^\infty e^{f(y)}\,dy\sim (e\pi)^{1/2}\,d^{1/4}\,
e^{-2\sqrt{d}},
\] 
as claimed.\qed

\textbf{Proof of Theorem \ref{clust}.}
Black vertices have the same local clustering coefficient in both
versions. Since the proportion of red vertices tends to be negligible
as $n\to\infty$, the limit of the average clustering coefficient is also
the same in both versions.
The global clustering coefficient of version 2 is identically equal to
$1$. In its defining fraction the numerator and the denominator are
proportional to $n$. When turning to version 1 the denominator have to
be increased by the number of triplets containing at least one red
edge. Such a triplet must have a red central vertex and at least one
more red vertex. Hence the increment of the denominator cannot exceed
$M_nZ_n^2$, where $M_n$ denotes the maximal degree, and $Z_n$ the
number of red vertices. In the proof of Proposition \ref{sharp} we
have shown that $M_n=O(n^{1/2+\varepsilon})$ and 
$Z_n=O(\log^{1+\varepsilon}n)$, thus the increment of the 
denominator is asymptotically negligible with respect to $n$. Hence
the global clustering   coefficient of version 1 must converge to
$1$. \qed

\end{document}